\documentclass[10pt]{article}

\usepackage{graphicx}
\usepackage{epsfig}

\usepackage{amsfonts,amsmath,amssymb}
\newcommand{\eps}{\varepsilon}
\newcommand{\cov}{\mathop{\mathsf{cov}}}

\newcommand{\E}{\mathsf{E}}
\newtheorem{lemma}{Lemma}

\newtheorem{theorem}{Theorem}

\newtheorem{remark}{Remark}
\newcommand{\Pp}{\mathsf{P}}
\newcommand{\R}{\mathbb{R}}
\newcommand{\Tb}{\mathbb{T}}
\newcommand{\Xb}{\mathbb{X}}
\newcommand{\N}{\mathbb{N}}

\newcommand{\Z}{\mathbb{Z}}

\newcommand{\ONE}{{\bf 1}}

\newcommand{\bpf}[1][Proof]{{\noindent {\sc #1: }}}
\newcommand{\epf}{{{\hspace{4 ex} $\Box$ \smallskip}}}

\author{Yuri Bakhtin\thanks{School of Mathematics, Georgia Tech, Atlanta GA, 30332-0160;  email:bakhtin@math.gatech.edu, 404-894-9235 (office phone), 404-894-4409(fax)}}

\title{Thermodynamic Limit for Large Random Trees.}
\begin{document}
\maketitle
\abstract{
We consider Gibbs distributions on finite random plane trees with bounded branching. We show that as the order of the tree grows to infinity, the distribution of any finite neighborhood of the root of the tree converges to a limit. We compute the limiting distribution explicitly and study its properties.
We introduce an infinite random tree consistent with these limiting distributions and show that it satisfies
a certain form of the Markov property. We also study the growth of this tree and prove several limit theorems including
a diffusion approximation.
}
\section{Introduction}
Various kinds of random trees have been studied in the literature. In this note we consider simply
generated random (plane rooted) trees also known as branching processes conditioned on the total population (CBP), see~\cite{Aldous-II:MR1166406}. Our
initial motivation was a study of the secondary structure statistics for large RNA molecules, see~\cite{Bakhtin-Heitsch-1:MR2415118}
and \cite{Bakhtin-Heitsch-2}. The secondary RNA structures can be encoded via plane rooted trees and studied with
the help of energy models. In~\cite{Bakhtin-Heitsch-1:MR2415118}
and \cite{Bakhtin-Heitsch-2}, it is demonstrated that the naive energy minimization approach to the prediction of typical secondary
structure features fails to explain the presence of high degree branchings. However, using the 
language of statistical mechanics and working with Gibbs ensembles on trees, we were able to include the entropy correction
and recover the typical RNA branching type. 
These results are concerned only with the rough information related to the branching statistics, but 
in this paper, we suggest a new viewpoint that helps to obtain some insights into the geometry of large random trees.

The model we work with follows the classical
Boltzmann--Gibbs postulate stating that the probablity of a configuration $T$ is proportional to $e^{-\beta E(T)}$, where
$E(T)$ is the energy of $T$, and $\beta$ is the inverse temperature in appropriate units (see the complete description of our model in Section~\ref{sec:setting_and_first_results}). Gibbs distributions, especially their limiting behaviour under the limit of the size of the system tending to infinity (so called thermodynamic limit),
are central to statistical mechanics, see \cite{Sinai:MR691854} and \cite{Georgii:MR956646} for a modern mathematical introduction.

The first goal of this paper is to prove that as the order of the tree grows to infinity, the distribution
induced by the Gibbs measure converges to that of an infinite discrete tree that we explicitly describe in detail (Sections~\ref{sec:setting_and_first_results} to~\ref{sec:limiting_tree}). This thermodynamic limit
belongs to the category of discrete limits of CBP according to the terminology introduced in~\cite{Aldous-II:MR1166406}, and our result (as well as the limiting object) appears to be new. In particular, it does not involve any rerooting procedures like the one introduced
in \cite{Aldous-I:MR1085326}.
We prove the result above for the bounded branching (or out-degree) case, although it should hold true under less restrictive assumptions.

The limiting infinite discrete tree is a more sophisticated object than a classical Galton--Watson tree. In particular, it dies out with zero probability and the progenies of distinct vertices are not independent. However, it turns out that the limiting tree is Markov in a natural sense, and the
Markov transition probability is explicitly computed in Section~\ref{sec:limiting_tree}. In Section~\ref{sec:growth_of_levels} we
notice that the number of vertices at a given distance $n$ from the root also form a Markov chain if $n$ is understood as a time parameter. We prove that under linear scaling
this Markov chain satisfies a limit theorem with the limit given by a gamma distribution. In Section~\ref{sec:flt} we strengthen this
result and show that a functional limit theorem holds with weak convergence to a diffusion process on the positive semi-line with constant drift and diffusion proportional to the square root of the space coordinate. Since this process (under the name of local time for Bessel(3) process) also serves as a scaling limit of the ``height profile'' for CBP itself, see~\cite[Conjecture 7]{Aldous-II:MR1166406} and~\cite{Gittenberger:MR1662793}, we can say that the infinite Markov random tree that we construct belongs to the same universality class as the original CBP.

There are several natural and interesting problems arising in connection with our results. One is, obviously, strengthening them to give an alternative to~\cite{Gittenberger:MR1662793} proof of the scaling limit in Aldous's Conjecture~7. Another one is to use our approach to study finer details of the random tree rather than the height profile. Our heuristic 
computation (see Section~\ref{sec:SPDE}) shows that the limit can be described as a solution of an SPDE with respect to a Brownian sheet.

{\bf Acknowledgements.} The author is grateful to NSF for partial support of this research via CAREER award DMS-0742424. He also thanks the
referees for their useful comments.

\section{The setting and first results on thermodynamic limit}\label{sec:setting_and_first_results}

\begin{figure}
\centering
\epsfig{file=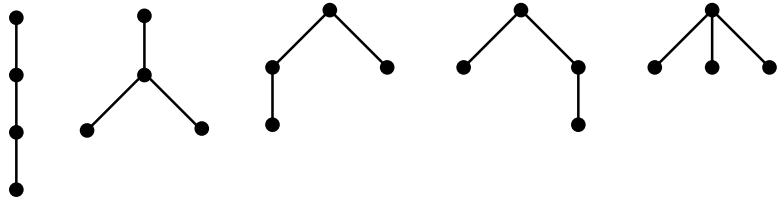,height=2.6cm}
\caption{Five different plane trees of order 4.}
\label{fig:4-trees}
\end{figure}

Let us recall that plane trees (or, ordered trees) are rooted trees
such that subtrees at any vertex are linearly ordered. In other words, two plane trees and are considered equal if there is a bijection between the
vertices of the two trees such that it preserves the parent --- child relation on the vertices and preserves the order of the child subtrees of any vertex.  Figure~\ref{fig:4-trees} shows all plane trees on $4$ vertices.

We fix $D\in\N$ and introduce $\Tb_N=\Tb_N(D)$, the set of all plane trees on~$N$ vertices such that the branching number (i.e.\ the number of children, or out-degree) of each vertex  does not exceed $D$. 
To introduce a Gibbs distribution on $\Tb_N$, we have to assign an energy value to each tree. We assume that an energy value 
$E_i\in\R$ is assigned to every $i\in\{0,\ldots,D\}$, and the energy of the tree $T$ is defined via
\[
 E(T)=\sum_{v\in V(T)} E_{\deg(v)}=\sum_{i=0}^D \chi_i(T)E_i,
\]
where $V(T)$ denotes the set of vertices of the tree $T$, $\deg(v)$ denotes the branching number of vertex~$v$, and $\chi_i(T)$ is the number of vertices of branching $i$ in $T$. Since the energy of an individual vertex depends only on its immediate neigborhood via the branching number,
one can say that this a model with nearest neighbor interaction.

Now we fix an inverse temperature parameter $\beta\in\R$ (usually, in statistical physics $\beta>0$, but our results apply to other values of $\beta$ as well) and define a probability measure $\mu_N$ on $\Tb_N$ by
\[
 \mu_N\{T\}=\frac{e^{-\beta E(T)}}{Z_N},
\]
 where the normalizing factor (partition function) is defined by
\[
Z_N=\sum_{T\in \Tb_N} e^{-\beta E(T)}.
\]
In particular, if $\beta=0$ or, equivalently, $E_i=0$ for all $i$, then $\mu_N$ is a uniform distribution on $\Tb_N$.

First, we are going to demonstrate that the above model admits a thermodynamic limit, i.e.\
the sequence of measures $(\mu_N)_{N\in\N}$ has a limit in a certain sense as $N\to\infty$. Secondly, we
study several curious properties of the limiting infinite random trees.

For each vertex $v$ of a tree $T\in\Tb_N$ its height $h(v)$ is defined as the distance to the root of $T$, i.e.\ the length of the shortest path connecting~$v$ to the root along the edges of $T$. The height of a finite tree is the maximum height of its vertices.

Let $n,N\in\N$. For any plane tree $T\in\Tb_N$, $\pi_{n,N} T$ denotes the neighborhood of the root of radius $n$, i.e.\ the subtree of $T$ spanned by all vertices with height not exceeding $n$. 

For any $n$ and sufficiently large $N$, the map $\pi_{n,N}$ pushes the measure~$\mu_N$ on~$\Tb_N$ forward to the measure $\mu_N\pi_{n,N}^{-1}$ on $S_n$, the set of all trees with height~$n$.

\begin{theorem}\label{th:main_convergence}
For each $n\in\N$, the measures $\mu_N\pi_{n,N}^{-1}$ on $S_n$ converge in total variation, as $N\to\infty$, to a measure $P_n$.
\end{theorem}
A proof of this theorem will be given in Section~\ref{sec:first_proofs}. At this point we prefer to introduce more definitions that will allow us to describe the limiting measures~$P_n$.

We define
\[
\Delta=\left\{p=(p_0,\ldots,p_D)\in [0,1]^{D+1}:\ \sum_{i=0}^Dp_i=1,\ \sum_{i=0}^Dip_i=1\right\},
\]
and let
\[
 J(p)= -H(p)+\beta E(p),\quad p\in\Delta.
\]
where 
\[
 H(p)=-\sum_{i=0}^{D}p_i\ln p_i
\]
is the entropy of the  probability vector $p\in\Delta$, and
\[
 E(p)=\sum_{i=0}^{D}p_iE_i
\]
is the associated energy.

The function $J$ is used to construct the rate function in the Large Deviation Principle for large plane trees, see \cite{Bakhtin-Heitsch-1:MR2415118},\cite{Bakhtin-Heitsch-2}.

It is strictly convex and its minimum value on $\Delta$ is attained at a unique point $p^*$. Using Lagrange's method, we
find that
\[
 \ln p_i^*+1+\beta E_i+\lambda_1+i\lambda_2=0,\quad i=0,1,\ldots,D,
\]
where $\lambda_1$ and $\lambda_2$ are the Lagrange multipliers. So we see that
\begin{equation}
 p^*_i=Ce^{-\beta E_i}\rho^i,\quad i=0,1,\ldots,D,
\label{eq:p*}
\end{equation}
where $C=e^{-1-\lambda_1}$, and $\rho=e^{-\lambda_2}$. 
In particular,
\begin{equation}
 p^*_i>0,\quad i=0,1,\ldots,D.
\label{eq:p*positive}
\end{equation}
Notice that $\rho$ can be characterized as a unique solution of
\begin{equation*}
\sum_{i=0}^De^{-\beta E_i}\rho^i=\sum_{i=0}^Die^{-\beta E_i}\rho^i,
\end{equation*}
and $C$ may be defined via
\begin{equation}
 \label{eq:C}\frac{1}{C}=\sum_{i=0}^De^{-\beta E_i}\rho^i=\sum_{i=0}^Die^{-\beta E_i}\rho^i.
\end{equation}
We denote $J^*=J(p^*)$ and $\sigma=e^{J^*}$.
For a tree $\tau\in S_n$, we introduce
\begin{equation}
 \bar E(\tau)=\sum_{\substack{v\in V(\tau)\\h(v)<n}}E_{\deg(v)}.
\label{eq:barE}
\end{equation}
Notice that the summation above excludes the highest level of the tree.

\begin{theorem}\label{th:limit_measure_evaluation} For any $n\in\N$, the limiting probability measure $P_n$ is given by
\begin{equation}
P_n\{\tau\}=Q_n k \rho^{k}\sigma^m e^{-\beta \bar E(\tau)}
\label{eq:P_n_up_to_Q_n}
\end{equation}
where the tree $\tau\in S_n$ is assumed to have $k$ vertices of height $n$ and $m$ vertices of height less than $n$. The constant $Q_n$
is a normalizing factor.
\end{theorem}
We give a proof of Theorems~\ref{th:main_convergence} and \ref{th:limit_measure_evaluation} in the next Section~\ref{sec:first_proofs}.
In Section~\ref{sec:Q} we compute the value of $Q_n$ explicitly. In Section~\ref{sec:limiting_tree} we shall see
that our convergence results may be interpreted as convergence to
an infinite random tree.

\begin{remark} Although Theorems~\ref{th:main_convergence} and \ref{th:limit_measure_evaluation} do not hold in full generality
for $D=\infty$, we expect that there is a large class of energy functions for which analogous results are true.
\end{remark}

\section{Proof of Theorems~\ref{th:main_convergence} and \ref{th:limit_measure_evaluation}}\label{sec:first_proofs}
For both theorems it is sufficient to check that for any $n$ and any two trees $\tau_1,\tau_2\in S_n$,
\begin{equation}
\label{eq:convergence_of_ratio}
 \lim_{N\to\infty}\frac{\mu_N\pi_{n,N}^{-1}\{\tau_1\}}{\mu_N\pi_{n,N}^{-1}\{\tau_2\}}=\frac{k_1e^{-\beta \bar E(\tau_1)}\rho^{k_1}\sigma^{m_1}}{
k_2e^{-\beta\bar  E(\tau_2)}\rho^{k_2}\sigma^{m_2}},
\end{equation}
where we assume that $\tau_1$ has $k_1$ vertices of height $n$, and $m_1$ vertices of height less than $n$;
$\tau_2$ has $k_2$ vertices of height $n$, and $m_2$ vertices of height less than~$n$. 

The energy of each tree $T$ with $\pi_{n,N} T=\tau_1$ is composed of contributions from the vertices of the tree $\tau_1$ of height less than $n$ (we
call this contribution $\bar E(\tau_1)$, see~\eqref{eq:barE}) and the contribution from the plane forest on $N-m_1$ vertices with $k_1$ connected components. The same applies to $\tau_2$.

Let us recall (see e.g.\ Theorem 5.3.10 in \cite{Stanley:MR1676282}) that the number of plane forests on $N$ vertices with $k$ components and $r_0,r_1,\ldots,r_D$ vertices with branching
numbers, respectively, $0,1,\ldots,D$ is
\[
 \frac{k}{N}\binom{N}{r_0,\ r_1,\ \ldots,\ r_D }
\]
if $r_0+\ldots+r_D=N$, $r_1+2r_2+\ldots+Dr_D=N-k$, and $0$ otherwise.

Therefore,
\begin{align}
 \frac{\mu_N\pi_{n,N}^{-1}\{\tau_1\}}{\mu_N\pi_{n,N}^{-1}\{\tau_2\}}&=
\frac{e^{-\beta \bar E(\tau_1)}\sum_{r\in\Delta(N,m_1,k_1)}\frac{k_1}{N-m_1}\binom{N-m_1}{r_0,\ r_1,\ \ldots,\ r_D }e^{-\beta E(r)}}{
e^{-\beta \bar E(\tau_2)}\sum_{r\in\Delta(N,m_2,k_2}\frac{k_2}{N-m_2}\binom{N-m_2}{r_0,\ r_1,\ \ldots,\ r_D }e^{-\beta
E(r)}}
\notag\\&=\frac{e^{-\beta\bar E(\tau_1)}I_1(N)}{e^{-\beta\bar E(\tau_2)}I_2(N)}.
\label{eq:ratio1}
\end{align}
Here
\begin{multline*}
\Delta(N,m,k)=\{r\in\Z_+^{D+1}:\ r_0+\ldots+r_D=N-m,\\ r_1+2r_2+\ldots+Dr_D=N-m-k\},
\end{multline*}
and $\Z_+=\N\cup\{0\}$.

Fix any $\eps>0$ and define 
\[
\Delta(N,m,k,\eps)=\left\{r\in \Delta(N,m,k):\ \left|\frac{r}{N-m}-p^*\right|<\eps\right\}.
\]

We claim that
\begin{equation}
\label{eq:Delta_equiv_Delta_eps}
 I_1(N)=I_1(N,\eps)(1+o(1)),\quad N\to\infty,
\end{equation}
where
\begin{equation*}
I_1(N,\eps)=\sum_{r\in\Delta(N,m_1,k_1,\eps)}\frac{k_1}{N-m_1}\binom{N-m_1}{r_0,\ r_1,\ \ldots,\ r_D }e^{-\beta E(r)}.
\end{equation*}
In fact, using Stirling's formula we see that if $r_i\ne 0$ for all $i=0,\ldots,D$,
\begin{align*}
 &\frac{k_1}{N-m_1}\binom{N-m_1}{r_0,\ r_1,\ \ldots,\ r_D }e^{-\beta E(r)}\\
&=\frac{k_1(N-m_1)^{N-m_1-\frac{1}{2}}e^{-\beta E(r)}e^{\frac{\theta_{N-m_1}}{12(N-m_1)}-\frac{\theta_{r_0}}{12r_0}-\ldots-\frac{\theta_{r_D}}{12 r_D}}}{
(2\pi)^{\frac{D}{2}} r_0^{r_0+\frac{1}{2}}\ldots r_D^{r_D+\frac{1}{2}}}\\ 
&=\frac{k_1 e^{-(N-m_1)J(\frac{r}{N-m_1})}}{((N-m_1)r_0\ldots r_D)^{\frac{1}{2}}}\cdot\frac{e^{\frac{\theta_{N-m_1}}{12(N-m_1)}-\frac{\theta_{r_0}}{12r_0}-\ldots-\frac{\theta_{r_D}}{12 r_D}}}{(2\pi)^{\frac{D}{2}}},
\end{align*}
with $0<\theta_j<1$ for all $j\in\N$. 
If $N$ is sufficiently large, there is a vector $r^*(N)\in \Delta(N,m_1,k_1,\eps)$ such that $|\frac{r^*(N)}{N-m_1}-p^*|<\eps/2$. Due to the strong convexity of $J$,
there is a number $\delta>0$ independent of $N$ such that
\[
 \min_{\Delta(N,m_1,k_1)\setminus\Delta(N,m_1,k_1,\eps)} J\left(\frac{r}{N-m_1}\right)>J\left(\frac{r^*(N)}{N-m_1}\right)+\delta,
\]
so that the contribution from each element of ~$\Delta(N,m_1,k_1)\setminus\Delta(N,m_1,k_1,\eps)$ is exponentially smaller than
that of $r^*(N)$ as $N\to\infty$.
The statement follows since the number of elements in~$\Delta(N,m_1,k_1)\setminus\Delta(N,m_1,k_1,\eps)$
is bounded by $N^{D+1}$. This argument can be easily extended to the case where $r_i=0$ for some $i$, which completes the proof of  
our claim~\eqref{eq:Delta_equiv_Delta_eps}.

Let us now define for $r\in \Delta(N,m_1,k_1,\eps)$,
\[
 b(r) = (r_0+(k_2-k_1),r_1-(k_2-k_1)-(m_2-m_1),r_2,r_3,\ldots,r_D ).
\]
Notice that for sufficiently small $\eps$ and sufficiently large~$N$, the image $\Delta'(N,\eps)$
of $\Delta(N,m_1,k_1,\eps)$ under $b$ is a subset of $\Delta(N,m_2,k_2)$. Moreover,
$b$ is invertible and, therefore, establishes a bijection between $\Delta(N,m_1,k_1,\eps)$ and $\Delta'(N,\eps)$.

Introducing
\[
I_2(N,\eps)=\sum_{r\in\Delta'(N,\eps)}\frac{k_2}{N-m_2}\binom{N-m_2}{r_0,\ r_1,\ \ldots,\ r_D }e^{-\beta E(r)},
\]
and using exactly the same reasoning as for $I_1$, we see that
\begin{equation}
 I_2(N)=I_2(N,\eps)(1+o(1)),\quad N\to\infty.
\label{eq:I_2(eps)}
\end{equation}
Equations \eqref{eq:ratio1},\eqref{eq:Delta_equiv_Delta_eps},\eqref{eq:I_2(eps)} imply now that
\begin{align}
  \frac{\mu_N\pi_{n,N}^{-1}\{\tau_1\}}{\mu_N\pi_{n,N}^{-1}\{\tau_2\}}&=\frac{e^{-\beta\bar E(\tau_1)}I_1(N,\eps)}{e^{-\beta \bar E(\tau_2)}I_2(N,\eps)}(1+o(1))\notag
\\
 &=\frac{k_1e^{-\beta\bar E(\tau_1)}}{k_2e^{-\beta\bar E(\tau_2)}}\cdot\frac{\sum_{r\in\Delta(N,m_1,k_1,\eps)}a_{1,r}}{\sum_{r\in\Delta(N,m_1,k_1,\eps)}a_{2,r}}(1+o(1)),\quad N\to\infty,
\label{eq:ratio_I_with_eps}
\end{align}
where
\[
 a_{1,r}=\binom{N-m_1}{r_0,\ r_1,\ \ldots,\ r_D }e^{-\beta E(r)},
\]
and
\[
 a_{2,r}=\binom{N-m_2}{r_0+(k_2-k_1),r_1-(k_2-k_1+m_2-m_1),r_2,\ldots,r_D}e^{-\beta E(b(r))}.
\]

Assuming that
$k_1\ge k_2$ and $m_1\ge m_2$ (all the other cases can be treated in the same way), we get
\[
\frac{a_{1,r}}{a_{2,r}}=\frac{((r_1-(k_2-k_1)-(m_2-m_1))\ldots(r_1+1)}{(N-m_2)\ldots(N-m_1+1))\cdot (r_0\ldots (r_0+(k_2-k_1)+1))}R,
\]
where
\[
 R=R(k_1,m_1,k_2,m_2)=e^{\beta(E_0-E_1)(k_2-k_1)-\beta E_1(m_2-m_1)}.
\]

Due to the definition of $\Delta(N,m,k,\eps)$,
\[
\frac{a_{1,r}}{a_{2,r}}\le \frac{((p_1^*+\eps)(N-m_1)-(k_2-k_1)-(m_2-m_1))^{-(k_2-k_1)-(m_2-m_1)}}{
(N-m_1)^{-(m_2-m_1)}((p_0^*-\eps)(N-m_1)+(k_2-k_1))^{-(k_2-k_1)}}R,
\]
so that
\begin{align}\notag
\limsup_{N\to\infty}\sup_{r\in\Delta(N,m,k,\eps)}\frac{a_{1,r}}{a_{2,r}}&\le (p_1^*+\eps)^{-(m_2-m_1)}
\left(\frac{p_1^*+\eps}{p_0^*-\eps}\right)^{-(k_2-k_1)}R\\
&\le \left(\frac{e^{-\beta E_1}}{p_1^*+\eps}\right)^{m_2-m_1}\left(\frac{(p_0^*-\eps)e^{\beta(E_0-E_1)}}{p_1^*+\eps}\right)^{k_2-k_1}
\label{eq:individual_terms1}
\end{align}
In the same way,
\begin{equation}
\label{eq:individual_terms2}
 \liminf_{N\to\infty}\inf_{r\in\Delta(N,m,k,\eps)}\frac{a_{1,r}}{a_{2,r}}\ge
\left(\frac{e^{-\beta E_1}}{p_1^*-\eps}\right)^{m_2-m_1}\left(\frac{(p_0^*+\eps)e^{\beta(E_0-E_1)}}{p_1^*-\eps}\right)^{k_2-k_1}
\end{equation}

Since the choice of $\eps$ is arbitrary, relations 
\eqref{eq:ratio_I_with_eps},\eqref{eq:individual_terms1}, and \eqref{eq:individual_terms2}
imply that
\begin{equation}
 \lim_{N\to\infty}\frac{\mu_N\pi_{n,N}^{-1}\{\tau_1\}}{\mu_N\pi_{n,N}^{-1}\{\tau_2\}}=\frac{k_1e^{-\beta\bar E(\tau_1)}}{k_2e^{-\beta\bar E(\tau_2)}}
\left(\frac{e^{-\beta E_1}}{p_1^*}\right)^{m_2-m_1}\left(\frac{p_0^*e^{\beta(E_0-E_1)}}{p_1^*}\right)^{k_2-k_1}.
\label{eq:limit_fraction}
\end{equation}

Using \eqref{eq:p*}, we see that
\begin{equation}
\frac{p_0^*e^{\beta(E_0-E_1)}}{p_1^*}=\frac{1}{\rho}.
\label{eq:p_0_p_1_rho}
\end{equation}
A direct computation based on \eqref{eq:p*} and~\eqref{eq:C} implies
\[
 H(p^*)=-\ln(C\rho)+\beta E(p^*).
\]
Therefore,
\begin{equation}
\label{eq:rhoC}
\frac{e^{-\beta E_1}}{p_1^*}=\frac{1}{C\rho}=e^{-J(p^*)}=\frac{1}{\sigma}.
\end{equation}
Now, \eqref{eq:convergence_of_ratio} is an immediate consequence of~\eqref{eq:limit_fraction},\eqref{eq:p_0_p_1_rho}, and \eqref{eq:rhoC}.
\epf

\section{Consistency and the precise value of~$Q_n$}\label{sec:Q}

We begin with the following consistency property:
\begin{theorem}\label{th:consistency}
The family of measures $(P_n)_{n\in\N}$ is consistent, i.e.\ for any $n$ and any $\tau\in S_n$
\[
P_n\{\tau\}=\sum_{\substack{\tau'\in S_{n+1}\\\pi_n^{n+1}\tau'=\tau}}P_{n+1}\{\tau'\},
\]
where $\pi_n^{n+1}$ denotes the projection map from  $S_{n+1}$ to $S_n$.
\end{theorem}
\bpf This theorem is a direct consequence of the limiting procedure in Theorem~\ref{th:main_convergence}. However, it is interesting to derive it from the specific form of $P_n$ provided by Theorem~\ref{th:limit_measure_evaluation}.

Let us assume that $\tau\in S_n$, and $\tau$ has $n$ vertices of height $k$ and $m$ of height less than $n$.
\begin{align*}
 &\sum_{\substack{\tau'\in S_{n+1}\\ \pi_n^{n+1}\tau=\tau}}P_{n+1}\{\tau'\}\\&=
 Q_{n+1}\sum_{i_1,\ldots,i_k=0}^De^{-\beta (\bar E(\tau)+E_{i_1}+\ldots+ E_{i_k} )}(i_1+\ldots+i_k)\rho^{i_1+\ldots+i_k}\sigma^{m+k}\\
&=Q_{n+1}e^{-\beta \bar E(\tau)}\sigma^{m+k}\sum_{i_1,\ldots,i_k=0}^D e^{-\beta(E_1+\ldots+E_k)}(i_1+\ldots+i_k)\rho^{i_1+\ldots+i_k}\\
&=Q_{n+1}e^{-\beta \bar E(\tau)}\sigma^{m+k}k\sum_{i_1=0}^D (i_1\rho^{i_1}e^{-\beta E_i})\sum_{i_2=0}^D
(\rho^{i_2}e^{-\beta E_{i_2}})\ldots \sum_{i_k=0}^D(\rho^{i_k}e^{-\beta E_{i_k}})\\
&=Q_{n+1}e^{-\beta \bar E(\tau)}\sigma^{m+k}k\frac{1}{C}\left(\frac{1}{C}\right)^{k-1}.
\end{align*}
In this calculation we denoted by $i_1,\ldots,i_k$ the branching numbers of the vertices of height $n$. We used the definition of $P_n$ in the first identity. The second identity is just a convenient rearrangement. The third one follows from the symmetry
in the factor $(i_1+\ldots+i_k)$. In the last identity we used~\eqref{eq:C} and the fact that $p^*\in \Delta$. Identity ~\eqref{eq:rhoC} implies
\begin{equation}
 \label{eq:1_over_C}
 \frac{1}{C}=\frac{\rho}{\sigma},
\end{equation}
so that
\begin{equation}
\label{eq:Q=sum_of_Q}
 \sum_{\substack{\tau'\in S_{n+1}\\ \pi_n^{n+1}\tau=\tau}}P_{n+1}\{\tau'\}=Q_{n+1}e^{-\beta \bar E(\tau)}\sigma^{m+k}k\sigma^{-k}\rho^k=\frac{Q_{n+1}}{Q_n}P_n\{\tau\}.
\end{equation}
Since this holds true for all $\tau\in S_n$, we can conclude that $Q_n=Q_{n+1}$ which completes the proof.\epf

Identity~\eqref{eq:Q=sum_of_Q} means that the constant $Q=Q_n$ in Theorem~\ref{th:limit_measure_evaluation} is the same for all $n$.
Choosing $n=1$ we can compute it using \eqref{eq:P_n_up_to_Q_n}:
\[
 1=Q\sum_{k=1}^D k e^{-\beta E_k}\rho^k\sigma^1=\frac{Q\sigma}{C}.
\]
A more precise version of Theorem~\ref{th:limit_measure_evaluation} easily follows:
\begin{theorem}\label{th:refined_limit_measure_evaluation} Let $C$ be defined by~\eqref{eq:C}. For each $n$, the limiting probability measure $P_n$ is given by
\[
P_n\{\tau\}= Ck e^{-\beta \bar E(\tau)}\rho^{k}\sigma^{m-1},
\]
where the tree $\tau\in S_n$ is assumed to have $k$ vertices of height $n$ and $m$ vertices of height less than $n$.
\end{theorem}

\section{The limiting random tree}\label{sec:limiting_tree}

Let $S_{\infty}$ be the set of infinite plane trees with branching number bounded by~$D$. Theorem~\ref{th:consistency} along with the classical Daniell---Kolmogorov Consistency theorem (see~\cite{Billingsley:MR1700749}) allows us to introduce a measure $P_{\infty}$ on $S_\infty$
consistent with measures $P_n$ for all $n$. Intuitively this is clear, but to make it precise we need to introduce a coding of plane trees. We have chosen one of several possible coding schemes.
 
 Let $T$ be a plane tree (finite or infinite) with branching bounded by $D$. Then~$T$ has a finite number $r_n\le D^n$ of vertices of any given height $n$. Let us say that all vertices of the same height $n$ form the $n$-th level of the tree. The vertices of $n$-th level are naturally ordered and can be enumerated by numbers from $1$ to $r_n$ (except for the case when there are no vertices at $n$-th level at all). Each of $r_n$ vertices of the $n$-th level has a parent at the level $n-1$. Denote the number received by the parent of $l$-th
vertex of the $n$-th level under the described enumeration by $g_{n,l}$. If $r_n<l\le D^n$ we set $g_{n,l}=0$. 
We also agree that for the root of the tree, i.e.,\ the first vertex in the zeroth level, $g_{0,1}=1$. 

Then for any $n\ge 0$ the $n$-th level  can be encoded by a vector \[g_n=(g_{n,1},\ldots,g_{n,D^n})\in\{0,1,\ldots,D^{n-1}\}^{D^n},\]
and the whole tree can be identified with the sequence of levels
\[
(g_1,g_2,\ldots)\in \Xb=\prod_{n=1}^{\infty}\{0,1,\ldots,D^{n-1}\}^{D^n},
\]
so that the space $\Tb$ of all plane trees (finite or infinite) with branching bounded by $D$ can be identified with a subset of $\Xb$. 

\begin{theorem}
 \label{th:measure_on_infinite_trees}
 There is a unique measure $P_{\infty}$ on $\Tb$ such that it is consistent with measures $P_n$:
\[
 P_{\infty}\pi_n^{-1}=P_n,
\]
where $\pi_n$ denotes the root's neighbourhood of height $n$ of a tree from $\Tb$.
This measure is concentrated on $S_\infty$.
\end{theorem}
\bpf The first statement follows from Theorem~\ref{th:consistency} and the Consistency theorem. The second statement
is a consequence of the fact that for each $n\in\N$, $P_n$ is concentrated on trees with positive number of vertices at $n$-th level.
\epf

The space $\Xb$ is compact in the product topology. Therefore, the convergence of finite-dimensional
distributions established in Theorem~\ref{th:main_convergence} and the classical Prokhorov theorem 
(see e.g.~\cite{Billingsley:MR1700749}) imply the following result:

\begin{theorem}\label{th:weak_convergence}
As $N\to\infty$, measures $P_N$ viewed as measures on $\Xb$ converge weakly to $P_\infty$ in the product topology. 
\end{theorem}
This statement shows that there is a limiting object for the random trees that we consider. This object
is an infinite random tree. For any $n\in\N$, the first $n$ levels of this random tree are distributed
according to $P_n$.

\begin{figure}
\centering
\epsfig{file=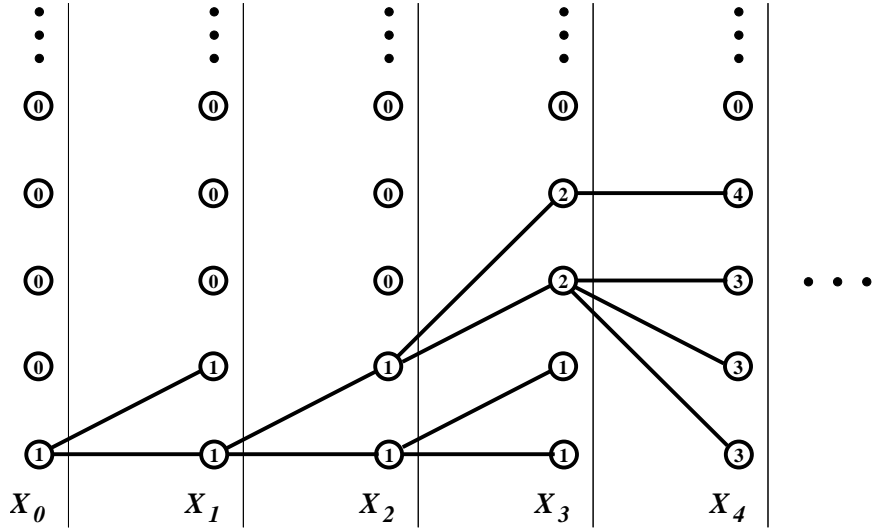,height=7cm}
\caption{A realization of the process $(X_n)$ and the associated random tree}
\label{fig:coding}
\end{figure}

Let us now embed the space $\Xb$ into $\bar \Xb=(\Z_+^{\N})^{\Z_+}$ filling up all the unused coordinates with zeros. The measure $P_\infty$ can be treated as a measure on 
$\bar\Xb$ thus generating a $\Z_+^{\N}$-valued process $(X_n)_{n=0}^\infty$ with discrete time. This process along with
the associated random tree is visualized
on Figure~\ref{fig:coding}. For any $n$, the map $X_n$ describes how the $n$-th level of the tree is built upon the $n-1$-th one.

For a level $g:\N\to\Z_+$, we denote by $|g|$ the number of non-zero entries in $g$ (i.e.\ the number of vertices at the level). 
For two levels~$g$ and~$g'$ we 
write $g\lhd g'$ if $\max_l g'_l\le |g|$.
If $g\lhd g'$ then we define 
\[
E(g,g')=\sum_{i=1}^{|g|} E_{\# \{j:\ g'_j=i\}},
\]
the energy induced by level $g'$ at its parent level $g$.

Theorem~\ref{th:refined_limit_measure_evaluation} immediately implies the following result:
\begin{theorem}\label{th:markov}
The process $(X_n)$ defined above is Markov with transition probability
\[
\Pp\{X_{n+1}=g'|\ X_{n}=g\}=\begin{cases}\frac{|g'|}{|g|}e^{-\beta E(g,g')}\rho^{|g'|-|g|}\sigma^{|g|},& g\lhd g',\\
                             0,&\mbox{\rm otherwise.}
                            \end{cases}
\]
\end{theorem}

\section{A limit theorem for the size of $n$-th level}\label{sec:growth_of_levels}

Let us introduce $Y_n=|X_n|$, the random number of vertices at $n$-th level. The following statement
is a direct consequence of Theorem~\ref{th:refined_limit_measure_evaluation} or Theorem~\ref{th:markov}:
\begin{theorem}
\label{th:markov-counting} The process $(Y_n)_{n=0}^{\infty}$ is Markov with transition
probability
\[
\Pp\{Y_{n+1}=k'|\ Y_{n}=k\}=\frac{k'}{k}\rho^{k'-k}\sigma^{k}\sum_{\substack{0\le i_1,\ldots,i_k\le D\\i_1+\ldots+i_k=k'}}e^{-\beta(E_{i_1}+\ldots+E_{i_k})}.
\]
\end{theorem}

The next theorem shows that in fact $Y_n$ grows
linearly in time. Let
\begin{equation}
 \mu=B_2-1,
\label{eq:mu}
\end{equation}
where
\begin{equation*}
 B_n=\sum_{i=0}^Di^np^*_i,\quad n\in\N.
\end{equation*}
Then $\mu>0$ being the variance of $p^*$, a nondegenerate distribution.

\begin{theorem}
\label{th:convergence_to_Gamma}
 \[
  \frac{Y_n}{n}\cdot\frac{2}{\mu}\stackrel{Law}{\to}\Gamma,
 \]
where $\Gamma$ is a random variable with density 
\[
p(t)=\begin{cases}te^{-t},&t\ge0\\ 0,&t<0\end{cases} 
\]
\end{theorem}

\bpf
Let us find the Laplace transform (generating function) of the distribution of $Y_n$:
\[
 L_n(s)=\E e^{sY_n},\quad s\le 0,
\]
(this definition differs from the traditional one by a sign change of the argument)
and prove that for any $x\le0$,
\begin{equation}
\label{eq:convergence_of_laplace}
\lim_{n\to\infty}L_n\left(\frac{x}{n}\right)=\frac{1}{\left(1-\frac{\mu x}{2}\right)^2}=:L_\infty(x),
\end{equation}
the r.h.s.\ being the Laplace transform of
\[
 p(t)=\frac{4t}{\mu^2}e^{-\frac{2t}{\mu}},
\]
the density of the r.v.\ $\frac{\mu}{2}\Gamma$.This will imply the desired result, 
see e.g.~\cite[Appendix 5]{Kallenberg:MR854102} for various statements on Laplace transforms.

Theorem~\ref{th:markov-counting} and~\eqref{eq:1_over_C} imply
\begin{align*}
&\E\left[e^{sY_{n+1}}|Y_n=k\right]=\sum_{k'}\frac{e^{sk'}k'}{k}\rho^{k'-k}\sigma^{k}\sum_{\substack{0\le i_1,\ldots,i_k\le D\\i_1+\ldots+i_k=k'}}e^{-\beta(E_{i_1}+\ldots+E_{i_k})}\\
&=\frac{\sigma^{k}}{k\rho^k}\sum_{0\le i_1,\ldots,i_k\le D}(i_1+\ldots+i_k)\rho^{i_1+\ldots+i_k}e^{s(i_1+\ldots+i_k)}e^{-\beta(E_{i_1}+\ldots+E_{i_k})}\\
&=w(s)v(s)^{k-1},
\end{align*}
where
\[
 v(s)=\sum_{i=0}^D p^*_ie^{si}=\sum_{i=0}^D C\rho^ie^{-\beta E_i}e^{si},
\]
and
\[
 w(s)=\sum_{i=0}^D ip^*_ie^{si}=\sum_{i=0}^D Ci\rho^ie^{-\beta E_i}e^{si}=v'(s).
\]
Therefore,
\begin{equation}
L_{n+1}(s)=\E w(s)v(s)^{Y_n-1}
          =\frac{w(s)}{v(s)}\E e^{\ln v(s) Y_n}
          =z(s)L_n(f(s)),\label{eq:L_iteration}
\end{equation}
where
\begin{equation}
\label{eq:taylor_for_f}
 f(s)=\ln v(s),
\end{equation}
and
\[
z(s)=\frac{w(s)}{v(s)}=f'(s).
\]

Both $z$ and $f$ are analytic functions. An elementary calculation shows that
\[
 f(s)=s+\frac{\mu}{2} s^2+r(s),
\]
and
\[
\ln z(s)=\mu s+q(s),
\]
where 
$\mu=w'(0)-1$ was introduced in~\eqref{eq:mu}
and 
\begin{equation}
\label{eq:remainders}
|r(s)|\le c|s|^3,\quad |q(s)|\le c|s|^2 
\end{equation}
for some $c>0$ and all $s\le 0$.

From now on, $x\le 0$ is fixed. Using \eqref{eq:L_iteration}
and  the obvious identity
\[
 L_0(s)=e^s,
\]
we can write
\[
L_n\left(\frac{x}{n}\right)=z\left(\frac{x}{n}\right)\cdot z\left(f\left(\frac{x}{n}\right)\right)\cdot z\left(f^2\left(\frac{x}{n}\right)\right)\cdot
\ldots \cdot z\left(f^{n-1}\left(\frac{x}{n}\right)\right)e^{f^{n}\left(\frac{x}{n}\right)},
\]
where 
\[
f^k(x)=\underbrace{f\circ f\ldots\circ f}_k(x),\quad k\ge 0,
\]
so that we have to study the numbers $(x_{n,k})_{n\in\N,k=0,\ldots,n}$ defined by
\[
 x_{n,k}={f}^{k}\left(\frac{x}{n}\right).
\]
We shall compare $(x_{n,k})_{n\in\N,k=0,\ldots,n}$ to $(y_{n,k})_{n\in\N,k=0,\ldots,n}$ defined by 
\[
 y_{n,k}=\frac{1}{\frac{n}{x}-\frac{\mu}{2}k}.
\]

For fixed $n$, both sequences $(x_{n,k})$ and $(y_{n,k})$ are negative and increasing in $k$. Therefore
\[
 |x_{n,k}|\le|x_{n,0}|=\frac{|x|}{n},
\]
and
\[
 |y_{n,k}|\le|y_{n,0}|=\frac{|x|}{n}.
\] 

Let us prove that for sufficiently large $n$ and any $k$ between~$0$ and~$n$, 
\begin{equation}
 |y_{n,k}-x_{n,k}|\le k\left(\frac{\mu^2}{4}+c\right)\left(\frac{|x|}{n}\right)^3.
\label{eq:y_close_to_x}
\end{equation}
This is certainly true for $k=0$. For the induction step, we write
\[
|y_{n,k}-x_{n,k}|\le|y_{n,k}-f(y_{n,k-1})|+|f(y_{n,k-1})-f(x_{n,k-1})|=I_1+I_2.
\]
A straightforward computation based on~\eqref{eq:taylor_for_f} shows that
\[
 |I_1|=\left|\frac{\mu^2}{4}y_{n,k-1}^2y_{n,k}+r(y_{n,k-1})\right|\le \left(\frac{\mu^2}{4}+c\right)\left(\frac{|x|}{n}\right)^3 
\]
Since $|f'(s)|\le 1$ for all sufficiently small $s$, we see that
\[
 |I_2|\le|y_{n,k-1}-x_{n,k-1}|.
\]
Combining these estimates we see that
\[
|y_{n,k}-x_{n,k}|\le \left(\frac{\mu^2}{4}+c\right)\left(\frac{|x|}{n}\right)^3+ |y_{n,k-1}-x_{n,k-1}|,
\]
and our claim~\eqref{eq:y_close_to_x} follows. It immediately implies that
\begin{equation}
 |y_{n,k}-x_{n,k}|\le \frac{C}{n^2}
\label{eq:y_close_to_x2}
\end{equation}
for some $K=K(x)$, sufficiently large $n$ and all $k$.
We can now write
\begin{align*}
 \ln L_n\left(\frac{x}{n}\right)&=\sum_{k=0}^{n-1}\ln z(x_{n,k})+x_{n,n}\\
&=\sum_{k=0}^{n-1}\mu x_{n,k}+\sum_{k=0}^{n-1}q(x_{n,k})+x_{n,n}\\
&=\sum_{k=0}^{n-1}\mu y_{n,k}+\sum_{k=0}^{n-1}\mu (x_{n,k}-y_{n,k})+\sum_{k=0}^{n-1}q(x_{n,k})+x_{n,n}\\
&=I_1+I_2+I_3+I_4.
\end{align*}
It is straightforward to see that $\lim_{n\to\infty}I_2+I_3+I_4=0$. The first term
\begin{align*}
I_1=\mu\sum_{k=0}^{n-1}\frac{1}{\frac{n}{x}-\frac{\mu}{2}k}
=\mu x\frac{1}{n}\sum_{k=0}^{n-1}\frac{1}{1-\frac{\mu x}{2}\frac{k}{n}}
\end{align*}
can be viewed as a Riemann integral sum, so that
\[
 \lim_{n\to\infty} \ln L_n\left(\frac{x}{n}\right)=\mu x \int_0^1\frac{du}{1-\frac{\mu x}{2}u}=-2\ln\left(1-\frac{\mu x}{2}\right),
\]
which immediately implies \eqref{eq:convergence_of_laplace}.
\epf

\section{A functional limit theorem}\label{sec:flt}
In this section we prove the following theorem on diffusion approximation for the process $Y$:
\begin{theorem} \label{th:FLT} Let
\[
 Z_n(t)=\frac{Y_{[nt]}}{n},\quad n\in\N, t\in\R_+.
\]
Then,  as $n\to\infty$, the distribution of $Z_n$ converges weakly in the Skorokhod topology in $D[0,\infty)$ to the
unique nonnegative weak solution $Z$ of the stochastic It\^o equation
\begin{align*}
 dZ(t)&=\mu dt+\sqrt{\mu Z(t)}dW(t),\\
 Z(0)&=0.
\end{align*}
\end{theorem}
\bpf Since the initial point $Z(0)=0$ is an ``entrance and non-exit'' singular point for the positive semi-axis
(see the classification of singular points in~\cite{Ito-Mckean:MR0345224} ), the
existence and uniqueness of a  nonnegative solution for all positive times is guaranteed.
Let us define
\[
b(x)\equiv \mu,\quad\text{and}\ a(x)=\mu\cdot\max\{x, 0\},\quad x\in\R,
\]
and extend the equation above to the negative semi-axis by
\begin{equation}
 dZ(t)=b(Z(t)) dt+\sqrt{a(Z(t))}dW(t).
\label{eq:extended equation}
\end{equation}
An obvious argument shows that there is no solution starting at $0$ and being negative for some $t>0$.
Therefore the weak existence and uniqueness in law hold for~\eqref{eq:extended equation}.
 According
to Section~5.4B of~\cite{Karatzas-Shreve:MR1640352}, this existence and uniqueness is equivalent to
the well-posedness of the martingale problem associated with $b$ and $a$.

We will use Theorem~4.1 from~\cite[Chapter 7]{Ethier-Kurtz:MR838085} on diffusion approximation. 
The coefficients $a,b$ were defined on the whole real line so as the theorem applies directly, with no modification. 
We proceed to check its conditions.

We must find processes $A_n$ and $B_n$ with the following properties:
\begin{enumerate}
\item  Trajectories of $A_n$ and $B_n$ are in $D[0,\infty)$.
\item  $A_n$ is nondecreasing.
\item $M_n=Z_n-B_n$ and $M_n^2-A_n$ are martingales with respect to the natural filtration generated by $Z_n,A_n,B_n$.
\item For every $T>0$ the following holds true:
\begin{align}
 \lim_{n\to\infty}\E\sup_{t\le T}|Z_n(t)-Z_n(t-)|^2&=0,
\\
 \lim_{n\to\infty}\E\sup_{t\le T}|A_n(t)-A_n(t-)|&=0,\label{eq:jumps_of_A}
\\
\lim_{n\to\infty}\E\sup_{t\le T}|B_n(t)-B_n(t-)|^2&=0,\label{eq:jumps_of_B}
\\
\sup_{t\le T}\left|B_n(t)-\int_0^tb(Z_n(s))ds\right|=\sup_{t\le T}\left|B_n(t)-\mu t\right|&\stackrel{\Pp}{\to}0,\quad n\to\infty,\label{eq:consistency_of_B_with_b}
\\
\sup_{t\le T}\left|A_n(t)-\int_0^ta(Z_n(s))ds\right|&\stackrel{\Pp}{\to}0,\quad n\to\infty.\label{eq:consistency_of_A_with_a}
\end{align}
\end{enumerate} 

We shall need the following lemma:
\begin{lemma}\label{lm:conditional_moments}
\begin{align*}
\E[Y_{j+1}|Y_j=k]=&\mu+k,\\
\E[Y_{j+1}^2|Y_j=k]=&B_3+3(k-1)B_2+(k-1)(k-2),\\
\E[Y_{j+1}^3|Y_j=k]=&B_4+4(k-1)B_3+6(k-1)(k-2)B_2+3(k-1)B_2^2\\&+(k-1)(k-2)(k-3),\\
\E[Y_{j+1}^4|Y_j=k]=&B_5+5(k-1)B_4+10(k-1)(k-2)B_3+10(k-1)B_3B_2
\\&+15(k-1)(k-2)B_2^2+10(k-1)(k-2)(k-3)B_2\\&+(k-1)(k-2)(k-3)(k-4). 
\end{align*}
\end{lemma}
\bpf For the first of these identities, we write
\begin{align*} 
\E[Y_{j+1}|\ Y_{j}=k]\notag
&=\frac{\sigma^{k}}{k\rho^k}\sum_{0\le i_1,\ldots,i_k\le D}(i_1+\ldots+i_k)^2\rho^{i_1+\ldots+i_k}e^{-\beta(E_{i_1}+\ldots+E_{i_k})}
\notag
\\
&=\frac{1}{k}\Biggl[k\left(\sum_{i_1=0}^D i_1^2 C\rho^{i_1}e^{-\beta E_{i_1}}\right)\left(\sum_{i_2=0}^D C\rho^{i_2}e^{-\beta E_{i_2}}\right)^{k-1}
\notag
\\
&+k(k-1)\left(\sum_{i_1=0}^D i_1 C\rho^{i_1}e^{-\beta E_{i_1}}\right)^2\left(\sum_{i_2=0}^D C\rho^{i_2}e^{-\beta E_{i_2}}\right)^{k-2}\Biggr]
\notag
\\
&=\frac{1}{k}(kB_2+k(k-1))=B_2+k-1
\\
&=\mu+k,
\end{align*}
where we used the symmetry of the terms $(i_1^2+\ldots i_k^2)$, $i_1i_2+i_1i_3+\ldots+i_{k-1}i_k$ 
and~\eqref{eq:1_over_C}. Next,
\begin{align*}
\E[Y_{j+1}^2|\ Y_j=k]
&=\frac{\sigma^{k}}{k\rho^k}\sum_{0\le i_1,\ldots,i_k\le D}(i_1+\ldots+i_{k})^3\rho^{i_1+\ldots+i_k}e^{-\beta(E_{i_1}+\ldots+E_{i_k})}
\\
&=\frac{1}{k}\bigl(k B_3+3k(k-1)B_2+k(k-1)(k-2)\bigr)
\\
&=B_3+3(k-1)B_2+(k-1)(k-2),
\end{align*}
and the other two identities in the statement of the lemma can be obtained in a similar way.
\epf

Returning to the proof of the functional limit theorem, let us find the coefficient $B_n(t)$ first. The process $Z_n$ is constant on any interval of the form $[j/n,(j+1)/n)$. Due to Lemma~\ref{lm:conditional_moments},
\begin{equation}
\label{eq:conditional_1st_moment}
 \E\left.\left[Z_n\left(t+\frac{1}{n}\right)\right|Z_n(t)\right]=Z_n(t)+\mu\frac{1}{n},
\end{equation}
so that we can set $B_n(t)=\mu[nt]/n$ to satisfy the martingale requirement on $M_n=Z_n-B_n$. Notice that with this
choice of $B_n$, relations~\eqref{eq:jumps_of_B} and~\eqref{eq:consistency_of_B_with_b} are easily seen to be satisfied.
Lemma~\ref{lm:conditional_moments} also implies
\begin{multline}
\label{eq:conditional_2nd_moment}
 \E\left.\left[Z^2_n\left(t+\frac{1}{n}\right)\right|Z_n(t)\right]
\\ 
=\frac{B_3+3(nZ_n(t)-1)B_2+(nZ_n(t)-1)(nZ_n(t)-2)}{n^2},
\end{multline}
so that for $t\in \frac{1}{n}\Z$,
\[
 \E\left.\left[M^2_n\left(t+\frac{1}{n}\right)-M^2_n(t)\right|Z_n(t)\right]=\frac{1}{n}\mu Z_n(t)+\frac{1}{n^2}(B_3-B_2^2-B_2+1).
\]
Therefore we can set
\[
 A_n(t)=\sum_{j:\frac{j}{n}\le t}\left(\frac{\mu}{n} Z_n\left(\frac{j}{n}\right)+\frac{[nt]}{n^2}(B_3-B_2^2-B_2+1)\right),
\]
to satisfy the martingale requirement on $M_n^2-A$. Notice that $A_n$ is nondecreasing since
\begin{align*}
\frac{1}{n}\mu Z_n(t)+\frac{1}{n^2}(B_3-B_2^2-B_2+1)&\ge \frac{1}{n^2}(B_2-1)+\frac{1}{n^2}(B_3-B_2^2-B_2+1)\\
&\ge\frac{1}{n^2}(B_3-B_2^2)\ge 0,
\end{align*} 
where the last inequality follows from the Cauchy---Schwartz inequality and $B_1=1$. So properties 1--3 are satisfied, 
and~\eqref{eq:consistency_of_A_with_a} follows from the definitions of $a$ and~$A$, and the convergence
\[
\lim_{n\to\infty} \sup_{t\le T}\sum_{j:\frac{j}{n}\le t}\frac{1}{n^2}(B_2-1)+\frac{1}{n^2}(B_3-B_2^2-B_2+1)=0.
\]

To prove \eqref{eq:jumps_of_A} we use the definition of $A$ to write
\[
\E\sup_{t\le T}|A_n(t)-A_n(t-)|\le \frac{\mu}{n}\E\sup_{\frac{j}{n}\le T} Z_n\left(\frac{j}{n}\right)+\frac{1}{n^2}(B_3-B_2^2-B_2+1)
\]
so it suffices to prove that $\E\sup_{\frac{j}{n}\le T} Z_n\left(\frac{j}{n}\right)$ is bounded.
The definition of~$M_n$, Lyapunov's inequality and Doob's maximal inequality for submartingales imply
that for some $c>0$:
\begin{align*}
\E\sup_{\frac{j}{n}\le T} Z_n\left(\frac{j}{n}\right)&\le \mu T+ \E\sup_{\frac{j}{n}\le T} \left|M_n\left(\frac{j}{n}\right)\right|
\\
&\le \mu T+ c\sqrt{\E M^2_n(T)}
\\
&\le \mu T+c\sqrt{ 2(\E Z^2_n(T)+\mu^2T^2)}.
\end{align*}

Lemma~\ref{lm:conditional_moments} implies that
$\E Z^2_n(T)$ has a limit, as $n\to\infty$, so that~\eqref{eq:jumps_of_A} is verified.

A lengthy but elementary calculation based on Lemma~\ref{lm:conditional_moments} shows that
\[
 \E[(Y_{j+1}-Y_j)^4|Y_n]\le c(Y_j^2+1)
\]
for some constant $c>0$, so that we can write
\begin{align*}
 \E\sup_{t\le T}(Z_n(t)-Z_n(t-))^2&\le \left[\E\sup_{t\le T}(Z_n(t)-Z_n(t-))^4\right]^{1/2}
\\ 
&\le \left[\frac{1}{n^4}\sum_{j:\frac{j}{n}\le T}\E(Y_{j+1}-Y_{j})^4\right]^{1/2}
\\
&
\le \left[\frac{c}{n^4}\sum_{j\le nT}\E(Y_{j}^2+1)\right]^{1/2}
\end{align*}
Since Lemma~\ref{lm:conditional_moments} implies that for some constant $c_1>0$,
\[
 \E(Y_j^2+1)\le c_1j^2,\quad j\in\N,
\]
we conclude that
\begin{align*}
 \E\sup_{t\le T}(Z_n(t)-Z_n(t-))^2
&\le
\sqrt{\frac{c }{n^4}\cdot n\cdot c_1 n^2T^2}\to 0,\quad n\to\infty,
\end{align*}
and the proof of the theorem is complete.
\epf

\section{Diffusion limit for finer structure of the random tree}\label{sec:SPDE}

In this section we present a non-rigorous and sketchy description for the diffusion limit of the infinite Markov random tree
itself rather then its width given by $Y_n$ at time $n$. Let us fix any time $n_0$ and divide all $Y_{n_0}$ vertices into $r$ nonempty disjoint groups. For any $n\ge n_0$ denote the progeny of $i$-th group at time $n$ by $V_{i,n}$.

We want to study the coevolution of $(V_{1,n},\ldots,V_{r,n})$. Though each~$V_{i,n}$ is not a Markov process, 
it is elementary to see that the whole vector is a homogeneous Markov process. We would like to compute the diffusion limit for this vector under an appropriate rescaling:
\[
 \frac{1}{n}(V_{1,[nt]},\ldots,V_{r,[nt]})
\]
We need to find the local drift and diffusion coefficients for the limiting process. Let $j_1+\ldots+j_r=k$.
Then computations similiar to Lemma~\ref{lm:conditional_moments} produce
\[
 \E\left[\frac{V_{1,m+1}}{n}-\frac{j_1}{n}\Bigr|\ \frac{1}{n}(V_{1,[nt]},\ldots,V_{r,[nt]})=\frac{1}{n}(j_1,\ldots,j_r)\right]=\mu\frac{j_1/n}{k/n}\frac{1}{n},
\]
so, by symmetry, the local limit drift is 
\[
b_i(v)=\mu \frac{v_i}{v_1+\ldots+v_r}.
\]
Similarly, the diagonal terms for local diffusion:
\begin{align*}
 \E&\left[\left(\frac{V_{1,m+1}}{n}-\frac{j_1}{n}\right)^2\Bigr|\ \frac{1}{n}(V_{1,[nt]},\ldots,V_{r,[nt]})=\frac{1}{n}(j_1,\ldots,j_r)\right]\\&=\mu\frac{j_1}{n}\frac{1+\frac{B_3-3B_2+2}{k}}{n},
\end{align*}
and
\[
 a_{ii}(v)=\mu v_i.
\]
For the off-diagonal terms a computation produces
\begin{align*}
\E[(V_{1,m+1}-j_1)(V_{2,m+1}-j_2)|\ V_{1,m}=j_1,\ldots,V_{r,m}=j_m]=0,
\end{align*}
so that
\[
 a_{ij}\equiv0,\quad i\ne j.
\]
So, the limiting equations are
\[
 dV_i(t)=\mu\frac{V_i(t)}{\sum_{j}V_j(t)}dt+\sqrt{\mu V_i(t)}\ONE_{\{V_i>0\}}dW_i(t).
\]
Let us introduce cumulative counts
\[
 U_j=V_1+\ldots+V_j.
\]
Then $Z(t)=U_r(t)$, and
\begin{multline}
\label{eq:dU}
 dU_j=\mu\frac{U_j}{U_r}dt+\sqrt{\mu U_1}\ONE_{\{U_1>0\}}dW_1+\ldots\\\ldots+\sqrt{\mu (U_j-U_{j-1})}\ONE_{\{U_j-U_{j-1}>0\}}dW_j
\end{multline}
Then, for each $0\le u_1\le\ldots\le u_r$, we can solve this equation with initial data
\[
(U_1(t_0),\ldots,U_r(t_0))=(u_1,\ldots,u_r),
\]
 which gives a random nondecreasing map
\begin{equation}
\label{eq:solution_map}
\Phi=\Phi_{t_0}:u\mapsto (U(t))_{t\ge t_0}.
\end{equation}
Here $u$ runs through the set $\{u_1,\ldots,u_r\}$. It is clear though that if we insert another point $u'$
between $0$ and $u^*=u_r$, then
solving the stochastic equation above for the modified set of initial points, we see that the 
the new solution map is a monotone extension of the old one. Adding points of a countable dense set one after another, we can extend
the solution map onto $u\in[0,u^*]$.
So, for each $u^*\ge0$ we are able to define a random monotone map $\Phi:[0,u^*]\to \R_+^{[t_0,\infty)}$.

Our last point is to represent these solution maps via stochastic integrals w.r.t.\ a  Brownian sheet
$(W(x,t))_{t,x\ge0}$, i.e. a continuous Gaussian random field with zero mean and 
\[\cov(W(x_1,t_1),W(x_2,t_2))=(x_1\wedge x_2)(t_1\wedge t_2),\quad x_1,x_2,t_1,t_2\ge0.\]

Equations \eqref{eq:dU} imply that $\Phi(u,t)$, $t\ge t_0$, $u\in[0,u^*]$ is equal in law to the monotonone (in $u$) solution of
the following SPDE:
\begin{align*}
d\Phi(u,t)=&\mu\frac{\Phi(u,t)}{\Phi(u^*,t)}dt+\int_{x\in\R}\ONE_{[0,\mu \Phi(u,t)]}W(dx\times dt),\\
\Phi(u,t_0)=&u,\quad u\in[0,u^*].
\end{align*}
A rigorous treatment of the limiting solution $\Phi$, and a precise convergence statement will appear elsewhere.
\bibliographystyle{alpha}
\bibliography{treelimit}
\end{document}